\documentclass[12pt,onecolumn]{article}
\usepackage[margin=1.1in]{geometry}
\usepackage{amsmath}
\usepackage{bm}	
\usepackage{graphicx}

\newtheorem{Theorem}{Theorem}

\newtheorem{Lemma}{Lemma}
\newtheorem{Proposition}{Proposition}

\usepackage{graphicx}
\usepackage{multirow}
\usepackage{mdwlist}

\usepackage{threeparttable}

\usepackage{amsfonts}
\usepackage{amsmath}
\usepackage{hyperref}
\usepackage{bbm}
\usepackage{subfigure}

\usepackage{enumerate}
\usepackage{multicol}
\usepackage{stackrel}

\usepackage{amssymb}
\usepackage[utf8]{inputenc}
\usepackage{hyperref}

\begin{document}
\title{\LARGE \bf Behavior and Management of Stochastic Multiple-Origin-Destination Traffic Flows Sharing a Common Link}

\author{Li Jin and Yining Wen
\thanks{This work was supported in part by NYU Tandon School of Engineering and C2SMART Department of Transportation Center. The authors appreciate discussion with Profs. Saurabh Amin and Dengfeng Sun.}
\thanks{L. Jin is with the Department of Civil and Urban Engineering and Y. Wen is with the Department of Mechanical and Aerospace Engineering, New York University Tandon School of Engineering, Brooklyn, NY, USA, emails: {lijin@nyu.edu, yw3997@nyu.edu}.}%
}
\newcommand*{\QEDA}{\hfill\ensuremath{\blacksquare}}%

\maketitle

\begin{abstract}
In transportation systems (e.g. highways, railways, airports), traffic flows with distinct origin-destination pairs usually share common facilities and interact extensively.
Such interaction is typically stochastic due to natural fluctuations in the traffic flows.
In this paper, we study the interaction between multiple traffic flows and propose intuitive but provably efficient control algorithms based on modern sensing and actuating capabilities.
We decompose the problem into two sub-problems: the impact of a merging junction and the impact of a diverging junction.
We use a fluid model to show that (i) appropriate choice of priority at the merging junction is decisive for stability of the upstream queues and (ii) discharging priority at the diverging junction does not affect stability.
We also illustrate the insights of our analysis via an example of management of multi-class traffic flows with platooning.
\end{abstract}


\section{Introduction}

In transportation systems such as roads \cite{kurzhanskiy2010active,osorio2013simulation,coogan2015compartmental,jin2018throughput} and airspace \cite{bertsimas1998air,sun2008multicommodity,chen2017stochastic}, traffic flows with distinct origin-destination pairs usually share common facilities (e.g. road section and airspace sector) to optimize system-wide efficiency and utilization of infrastructure. Consequently, multiple traffic flows interact extensively in the common link, and such interaction can propagate to upstream links.

\begin{figure}[hbt]
\centering
\includegraphics[width=0.45\textwidth]{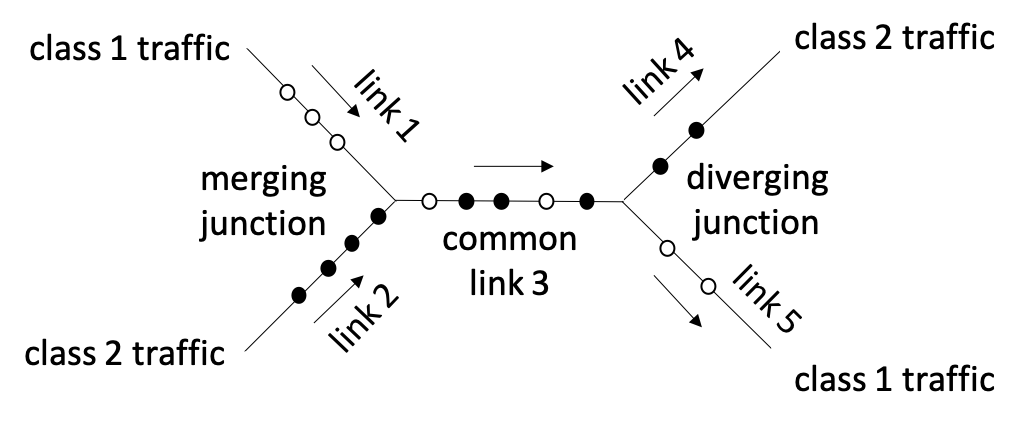}
\caption{Two traffic flows with distinct origin-destination pairs sharing a common highway section.}
\label{road}
\end{figure}

Consider the typical setting in Fig.~\ref{road}. Two classes of traffic ``compete'' for getting discharged from the common link 3, which can lead congestion in link 3. 
Limited capacities of links 4 and 5 can also contribute to this congestion.
Congestion in link 3 may further block traffic from the upstream links.

If both inflows at the source nodes and capacities of the links are constant, then no congestion should arise as long as the inflows are less than the capacities.
However, in reality, congestion is prone to occur due to fluctuations in inflows.
For example, inflows to a highway depends on traffic condition on upstream arterial roads as well as demand-disrupting events (concert, sports, etc.)
Air traffic flow is heavily influenced by weather.
Furthermore, such fluctuation is typically stochastic and is best modeled probabilistically.
However, how to manage traffic flows in such scenarios, especially under stochastic inflows, has not been well understood.

In this paper, we study the behavior of multiple traffic flows sharing a common link and propose intuitive but provably efficient management strategies that ensure bounded queuing delay and maximal throughput.
We consider the setting where both links 1 and 2 are subject to Markovian inflows: the inflow to each link switches between two values according to a Markov chain.
We assume that the inflows to these two links are statistically independent.
The traffic flows have their respective fixed routes, which overlap on link 3.
Link 3, the common link, has a finite storage space; once the traffic queue in link 3 attains the storage space, the flows out of links 1 and 2 will be reduced due to spillback.
The limited storage space of link 3 is shared by traffic from links 1 and 2 according to pre-specified priorities.
The multiclass traffic flow is discharged from link 3 according to a proportional rule: the discharge rate of a traffic class is proportional to the fraction of traffic of this class in the current queue.
Such discharging rule may also cause spillback from links 4 or 5 to link 3.

In this setting, the major decision variable for traffic flow management is the priorities according to which the limited capacity of the common link is shared between two traffic classes. In road traffic, this involves signal control (urban streets \cite{osorio2013simulation}) and ramp metering (highways \cite{kurzhanskiy2010active}). In air transportation, this involves ground and/or airborne holding \cite{sun2008multicommodity,zhou2019resilient}. 
Such control actions are typically costly and must be designed based on rigorous and systematic justification.

The main contributions of this paper are a set of results that help a system operator determine the priorities at the merging junction based on operational parameters (demands and capacities) to ensure guarantees of key performance metrics, viz. queuing delay and throughput.
Specifically, we argue via Theorem~\ref{thm_1} that there exist a non-empty set of priorities ensuring bounded traffic queues at the merging junction and maximal throughput if and only if the average inflow of each traffic class is less than the capacity of each link on its route.
We further argue via Theorem~\ref{thm_2} that the discharging rule and the possible spillback at the diverging junction does not affect stability of the system.
In addition, we explicitly provide a set of priorities that stabilize the system.
We expect our results to be directly relevant for road traffic management \cite{jin2018stability} and the discrete-state extension thereof to be relevant for air traffic management \cite{chen2017stochastic}.

Our study is based on a fluid model.
Fluid models are commonly used for highway bottlenecks \cite[Ch. 2]{newell13}.
Its discretized version is also common for air traffic management \cite{bertsimas1998air}.
We are aware that queuing models (e.g. M/M/1) are also widely used in transportation studies \cite{osorio2013simulation,baykal2009modeling}.
However, queuing models focus on the delay due to random headways between vehicles rather than the congestion due to demand fluctuation; therefore, fluid model fits our purpose better.
In fact, we view our fluid model as a reduced-order abstraction for queuing model:
it is well known that stability of queuing models is closely related to their fluid counterparts \cite{dai1995stability,bertsimas1996stability}.
Hence, our analysis in itself contributes to the literature on stochastic fluid models, which has mainly focused on controlling single/parallel links \cite{jin2018stability,kulkarni97,cassandras02} or characterizing steady-state distribution of queue sizes \cite{mitra1988stochastic,kroese2001joint}, rather than quantifying spillback-induced delay and throughput loss.

The rest of this paper is organized as follows.
In Section~\ref{sec_merge}, we isolate the merging junction from the network and study its behavior.
In Section~\ref{sec_system}, we add the diverging junction into our analysis and obtain results for the merge-diverge system.
In Section~\ref{sec_simulate}, we present a numerical example illustrating the main results that we derived.
In Section~\ref{sec_conclusion}, we summarize the main conclusions and propose several directions for future work. 
\section{Analysis of merging junction}
\label{sec_merge}

\begin{figure}[hbt]
\centering
\includegraphics[width=0.3\textwidth]{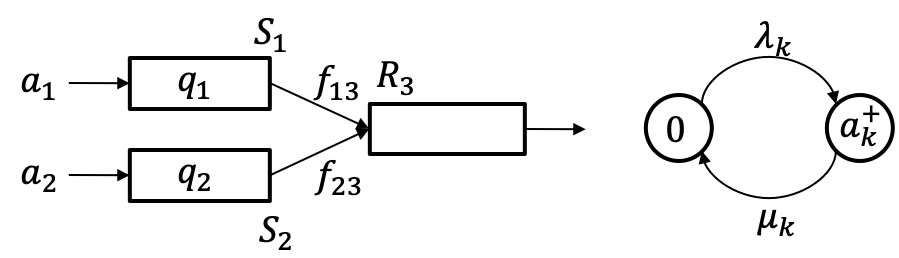}
\caption{Merging junction (left); inflow to link $k\in\{1,2\}$ evolves according to a Markov chain (right).}
\label{y}
\end{figure}

In this section, we study the behavior of a single merging junction (see Fig.~\ref{y}). This is an important component in the merge-diverge system and turns out to play a decisive role in terms of stability and throughput analysis.

\subsection{Model and main result}
\begin{figure}[hbt]
\centering
\includegraphics[width=0.15\textwidth]{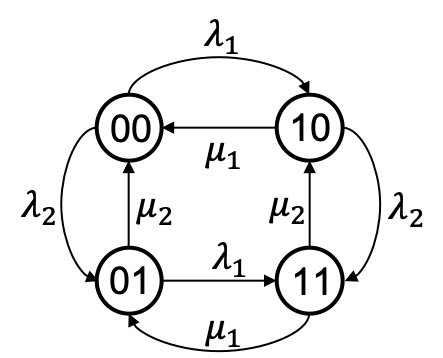}
\caption{Markov chain for $\{A(t)=[A_1(t)\ A_2(t)]^T;t>0\}$.}
\label{chain}
\end{figure}

The merging junction consists of three links. Traffic flows out of the upstream links 1 and 2 join and enter the downstream link 3.
The \emph{inflows} to links 1 and 2 are specified as follows. Let $A_k(t)$ be the inflow to link $k$ at time $t$. Then, $A_k(t)$ is a two-state Markov process with state space $\{0,a_k^+\}$ and transition rates $\lambda_k$ and $\mu_k$, as illustrated in Fig.~\ref{y}.
Thus, the mean inflows are given by
\begin{align*}
\bar a_k=\frac{\lambda_k}{\lambda_k+\mu_k}a_k^+
\quad k=1,2.
\end{align*}
Note that our results can be extended to the case where $A_k(t)$ is a two-state Markov process with state space $\{a_k^-,a_k^+\}$ for some $a_k^->0$.
We assume that $\{A_1(t);t>0\}$ and $\{A_2(t);t>0\}$ are independent processes. Consequently, we can use a four-state Markov chain to describe the evolution of the inflows. The state of the Markov chain is $a\in\mathcal A=\{0,a_1^+\}\times\{0,a_2^+\}$. Fig.~\ref{chain} uses a shorthand notation where ``00'' means $a=[0\ 0]^T$ and ``10'' means $a=[a_1^+\ 0]^T$.
We also use the unified notation $\{\nu_{ij};i,j\in\mathcal A\}$ to denote the transition rates; e.g. $\nu_{00,10}=\lambda_1$.

The \emph{flows} $f_{13},f_{23}$ between the links are determined by the \emph{sending flows} offered by links 1 and 2 as well as the \emph{receiving flow} allowed by link 3. Specifically, let $q_k\in[0,\infty)$ be the \emph{queue length} in link $k$ and $\mathcal Q=[0,\infty)^2$ be the set of queue lengths. Then, the sending flow out of link $k$ is given by
\begin{align*}
s_k(q_k,a_k)=\begin{cases}
a_k & q_k=0,\\
F_k & q_k>0,
\end{cases}
\quad k=1,2,
\end{align*}
where $F_k$ is the \emph{capacity} of link $k$.
The receiving flow of link 3 is given by
\begin{align*}
r_3(q_3)=\begin{cases}
R_3 & q_3<\theta,\\
F_3 & q_3=\theta,
\end{cases}
\end{align*}
where $R_3$ is the maximal receiving flow of link 3.
In this section, we focus on the merging junction and assume that link 3 is not constrained downstream and $Q_3(t)=0$ for all $t\ge0$; we will relax this assumption in the next section.
Thus, the between-link flows are given by
\begin{subequations}
\begin{align}
&f_{13}(a,q)=\min\{s_1(a_1,q_1),R_3\mathbb I_{q_2=0}+\phi_1 R_3\mathbb I_{q_2>0}\},\label{eq_f1}\\
&f_{23}(a,q)=\min\{s_2(a_2,q_2),R_3\mathbb I_{q_1=0}+\phi_2 R_3\mathbb I_{q_1>0}\}\label{eq_f2}
\end{align}
\end{subequations}
where $\phi_k\in[0,1]$ is the \emph{priority} of link $k=1,2$ and $R_3=r_3(0)$ is the receiving flow of link 3 for $q_3=0$; the \emph{priority vector} $\phi=[\phi_1\ \phi_2]^T$ must satisfy
\begin{align}
\phi\ge0,\quad|\phi|=1.
\label{eq_phi}
\end{align}
In practice, the priority vector specify how the limited capacity of the shared link is distributed over the upstream links. A typical mechanism for implementing such capacity allocation is traffic signal control, i.e. intersection control on urban streets and ramp metering on highways. In air transportation, this is done by air traffic management instructions.

The \emph{state} of the merging junction is $(a,q)\in\mathcal A\times\mathcal Q$. The evolution of $A(t)$ is fully specified by the Markov chain in Fig.~\ref{chain}. For a given initial condition $Q(0)\in\mathcal Q$, the evolution of $Q(t)$ is governed by
\begin{align*}
\frac{d}{dt}Q_k(t)=A_k(t)-f_{k3}(A(t),Q(t))
\quad k=1,2,\
t>0.
\end{align*}
The process $\{(A(t),Q(t));t>0\}$ is actually a piecewise-deterministic Markov process \cite{davis84}.

We say that the the merging junction is \emph{stable} if the total queue size is bounded, i.e. if there exists $Z<\infty$ such that for each initial condition
\begin{align}
\limsup_{t\to\infty}\frac1t\int_{s=0}^t\mathsf E[Q_1(s)+Q_2(s)]ds\le Z.
\label{eq_stable}
\end{align}
This definition of stability is motivated by \cite{dai95}.
Note that stability of the merging junction depends on the inflow, the capacities, the maximal receiving flow, and the priority vector.

The main result of this section is a criterion for existence of priority vectors that stabilize the queues:

\begin{Theorem}
\label{thm_1}
Consider a merging junction and let $[\phi_1\ \phi_2]^T\in[0,1]^2$ be the priority vector satisfying \eqref{eq_phi}. Then, there exists a non-empty set of priority vectors that stabilize the merging junction if and only if
\begin{align}
\bar a_1<F_1,\
\bar a_2<F_2,\
\bar a_1+\bar a_2<R_3.
\label{eq_nominal}
\end{align}
Furthermore, when \eqref{eq_nominal} holds, if furthermore
\begin{align}
\frac{\bar a_1}{F_1}+\frac{\bar a_2}{F_2}<1
\label{eq_sufficient1}
\end{align}
holds, then every $\phi\in\Phi$ is stabilizing; otherwise, a set of stabilizing priority vectors is given by
\begin{align}
\Phi_1=\{\phi\in\Phi:\phi_1>\bar a_1/R_3,\ \phi_2>\bar a_2/R_3\}
\label{eq_Phi1}
\end{align}
and a set of destabilizing priority vectors is given by the complement of the set
\begin{align*}
\Phi_0=\Bigg\{&\phi\in\Phi:\frac{\bar a_1}{F_1}+\frac{\bar a_2}{F_2}+\left(1-\frac{\phi_1 R_3}{F_1}-\frac{\phi_2R_3}{F_2}\right)\\
&\times\min\left\{\frac{\bar a_1}{\phi_1 R_3},\frac{\bar a_2}{\phi_2R_3}\right\}\le1\Bigg\}.
\end{align*}
\end{Theorem}

The above theorem essentially states that there exist stabilizing priority vectors if and only if the average inflows are less than the respective capacities.
Furthermore, we provide criteria for the stability of particular priority vectors.
Note that the set $\Phi_1$ (resp. $\Phi_0$) is derived from a sufficient (resp. necessary) condition for stability; there may exist a gap between $\Phi_1$ and $\Phi_0$.
For priority vectors in the gap, our results do not provide a conclusive answer regarding stability; see Section~\ref{sec_simulate} for a numerical example with a visualization of this gap.

\subsection{Proof of Theorem~\ref{thm_1}}
This subsection is devoted to a series of results leading to Theorem~\ref{thm_1}.
First, we derive a necessary condition for stability of the merge:

\begin{Proposition}
\label{prp_necessary}
Consider the merging junction and let $\phi_1\in[0,1]$ be the priority of link 1. If the traffic queues upstream to the merging junction are stable, then either
\begin{align}
\frac{\bar a_1}{F_1}+\frac{\bar a_2}{F_2}\le1
\label{eq_necessary1}
\end{align}
or
\begin{align}
&\frac{\bar a_1}{F_1}+\frac{\bar a_2}{F_2}+\left(1-\frac{\phi_1 R}{F_1}-\frac{\phi_2R}{F_2}\right)\min\left\{\frac{\bar a_1}{\phi_1 R},\frac{\bar a_2}{\phi_2R}\right\}\le1
\label{eq_necessary2}
\end{align}
holds, where $1/\phi_k:=\infty$ for $\phi_k=0$, $k=1,2$.
\end{Proposition}

\emph{Proof:}
Apparently, for each initial condition, the fluid queuing process $\{Q(t);t>0\}$ can always visit the state $q_1=0,q_2=0$ within finite time and with positive probability. Hence, the fluid queuing process is ergodic \cite[Theorem 2.11]{cloez2015exponential}. Hence, there exist constants $\mathsf p_{01}$, $\mathsf p_{10}$, and $\mathsf p_{11}$ such that for any initial condition
\begin{align*}
&\mathsf p_{00}=\lim_{t\to\infty}\frac1t\int_{s=0}^t\mathbb I_{Q_1(s)=0,Q_2(s)=0}ds\quad a.s.\\
&\mathsf p_{01}=\lim_{t\to\infty}\frac1t\int_{s=0}^t\mathbb I_{Q_1(s)=0,Q_2(s)>0}ds\quad a.s.\\
&\mathsf p_{10}=\lim_{t\to\infty}\frac1t\int_{s=0}^t\mathbb I_{Q_1(s)>0,Q_2(s)=0}ds\quad a.s.\\
&\mathsf p_{11}=\lim_{t\to\infty}\frac1t\int_{s=0}^t\mathbb I_{Q_1(s)>0,Q_2(s)>0}ds\quad a.s.
\end{align*}
where $\mathbb I$ is the indicator variable.
If the upstream queues are stable, then
\begin{subequations}
\begin{align}
&\bar a_1=\mathsf p_{10}F_1+\mathsf p_{11}\phi_1 R\label{eq_a1}\\
&\bar a_2=\mathsf p_{01}F_2+\mathsf p_{11}\phi_2 R
\end{align}
\end{subequations}
Also note that
\begin{subequations}
\begin{align}
&\mathsf p_{01}\ge0,\
\mathsf p_{10}\ge0,\
\mathsf p_{11}\ge0,\\
&\mathsf p_{01}+\mathsf p_{10}+\mathsf p_{11}\le1.\label{eq_ppp}
\end{align}
\end{subequations}
Then, one can obtain from \eqref{eq_a1}--\eqref{eq_ppp} that either \eqref{eq_necessary1} or \eqref{eq_necessary2} holds.
\hfill$\blacksquare$

Second, we derive a sufficient condition for stability of the merge:

\begin{Proposition}
\label{prp_sufficient}
Consider a merging junction and let $\phi_1\in[0,1]$ be the priority of link 1. The traffic queues upstream to the merging junction are stable if either (i)
\eqref{eq_nominal}--\eqref{eq_sufficient1} hold
or (ii) the following inequalities
\begin{subequations}
\begin{align}
&\bar a_1-\min\{F_1,\phi_1 R\}<0
\label{eq_sufficient2a}\\
&\bar a_2-\min\{F_2,\phi_2 R\}<0.
\label{eq_sufficient2b}
\end{align}
\end{subequations}
hold.
\end{Proposition}

\emph{Proof:}
We only prove the case where \eqref{eq_sufficient1} holds; the case where \eqref{eq_sufficient2a}--\eqref{eq_sufficient2b} hold can be proved analogously. Suppose that \eqref{eq_sufficient1} holds. Consider the quadratic Lyapunov function
\begin{align*}
V_1(a,q):=&q^T
\left[\begin{array}{cc}
1 & \alpha\\
\alpha & \alpha^2
\end{array}\right]
q
+[\beta_i\ \alpha\beta_i]q,
\\
&\quad q\in[0,\infty)^2,\ i\in\{00,10,01,11\}
\end{align*}
where the parameters are given by
\begin{subequations}
\begin{align}
&\alpha:=\frac12\left(\frac{\bar a_1}{F_2-\bar a_2}+\frac{F_1-\bar a_1}{\bar a_2}\right),\\
&\beta_{00}: =1,\
\beta_{10}:=\frac{\bar a_1}{\lambda_1}+1,\label{eq_beta00}\\
&\beta_{01}:=\alpha\frac{\bar a_2}{\lambda_2}+1,\
\beta_{11}:=\frac{\bar a_1}{\lambda_1}+\alpha\frac{\bar a_2}{\lambda_2}+1.\label{eq_beta11}
\end{align}
\end{subequations}
The above parameters are guaranteed to be positive by \eqref{eq_sufficient1}.
Let $\mathcal L_1$ be the infinitesimal generator (see \cite{davis84} for definition) of the merge system. With $\alpha$ and $\beta_i$ as given above, one can show that for each $i$ and $q$
\begin{align*}
\mathcal L_1V_1(i,q)=\Big((\bar a_1-f_1(q))+\alpha(\bar a_2-f_2(q))\Big)(q_1+\alpha q_2).
\end{align*}
Then, there exist
\begin{align*}
&c:=\min\{F_1-\bar a_1-\alpha\bar a_2,\alpha F_2-\bar a_1-\alpha\bar a_2\}\min\{1,\alpha\}>0\\
&d:=\max_{q_1\le\tilde q_1,q_2\le\tilde q_2}\mathcal L_1V_1(i,q)<\infty
\end{align*}
such that
\begin{align*}
\mathcal L_1V_1(i,q)\le-c|q|+d
\quad\forall i,q
\end{align*}
which implies \eqref{eq_stable} according to the Foster-Lyapunov drift criteria \cite[Theorem 4.3]{meyn93}.
\hfill$\blacksquare$

Finally, we optimize the priority under various scenarios. The following result characterizes the priorities that leads to maximal throughput under the stability constraint given by Proposition~\ref{prp_sufficient}:

\emph{Proof of Theorem~\ref{thm_1}:}
The necessity of \eqref{eq_nominal} is apparent: if the queues are stable, then the average inflows must be less than the corresponding capacities.

To show the sufficiency of \eqref{eq_nominal}, note that \eqref{eq_nominal} implies that every $\phi\in\Phi_1$ verifies \eqref{eq_sufficient2a}--\eqref{eq_sufficient2b}. The sets $\Phi_1$ and $\Phi_1'$ naturally result from Propositions~\ref{prp_necessary} and \ref{prp_sufficient}.
\hfill$\blacksquare$ 
\section{Analysis of merge-diverge network}
\label{sec_system}

\begin{figure}[hbt]
\centering
\includegraphics[width=0.3\textwidth]{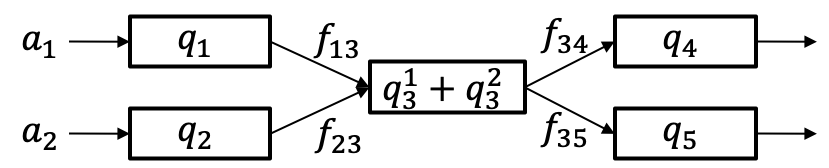}
\caption{Merge and diverge in series; two classes of traffic flows share a common link.}
\label{yy}
\end{figure}

In this section, we extend the results for a merging junction to the merge-diverge network in Fig.~\ref{yy}.

\subsection{Model and main result}
The dynamics is essentially the same as that of the merging junction; the main difference results from the interaction between two traffic classes at the diverge according to the \emph{discharging rule}, which is specified by $\psi$, the fraction of class 1 traffic in the sending flow $s_3$.
For the flow out of link 3, we consider a \emph{proportional} discharging rule $\psi:[0,\infty)^2\times\prod_{k=1}^2[0,F_k]\to[0,1]^2$ defined as
\begin{align}
&\psi(q_3^1,q_3^2,f_{13},f_{23})\nonumber\\
&=\begin{cases}
\left[\frac{q_3^1}{q_3^1+q_3^2}\ \frac{q_3^2}{q_3^1+q_3^2}\right]^T & q_3^1+q_3^2>0\\
\left[\frac{f_{13}}{f_{13}+f_{23}}\ \frac{f_{23}}{f_{13}+f_{23}}\right]^T & q_3^1+q_3^2=0,f_{13}+f_{23}>0\\
[1/2\ 1/2]^T & o.w.
\end{cases}
\label{eq_psi}
\end{align}

The sending and receiving flows are $s_k(q_k)$ and $r_k(q_k)$ respectively, $k=1,2,\ldots,5$. With a slight abuse of notation, we denote $q=[q_1\ q_2\ \cdots q_5]^T$. The inter-link flows are given by
\begin{subequations}
\begin{align}
&f_{13}(q)=\min\left\{s_1(q_1),\Big(r_3(q_3)-(1-\phi)s_2(q_2)\Big)_+\right\}
\label{eq_f13}\\
&f_{23}(q)=\min\left\{s_2(q_2),\Big(r_3(q_3)-\phi s_1(q_1)\Big)_+\right\}
\label{eq_f23}\\
&f_{34}(q)=\min\left\{\psi s_3(q_3),r_4(q_4),\frac{\psi}{1-\psi}r_5(q_5)\right\}
\label{eq_f34}\\
&f_{35}(q)=\min\left\{\psi s_3(q_3),r_5(q_5),\frac{1-\psi}{\psi}r_4(q_4)\right\}
\label{eq_f35}
\end{align}
\end{subequations}
In the above, \eqref{eq_f13}--\eqref{eq_f23} are direct extensions of \eqref{eq_f1}--\eqref{eq_f2} from $q_3=0$ to $q_3\ge0$. \eqref{eq_f34}--\eqref{eq_f35} essentially follow the same logic; the only difference is the impact of the discharging ratio $\psi$. That is, congestion in one traffic class (e.g. $q_5$) may undermine the flow in the other (e.g. $f_{34}$).
For ease of presentation, we assume that 
$$R_4<F_3,\ R_5<F_3,\ F_3<R_4+R_5.$$
The other cases can be analogously studied but are less interesting.

Since the queue in link 3 is upper-bounded and since links 4 and 5 are not constrainted downstream, the merge-diverge network is stable if \eqref{eq_stable} holds. However, in this setting the upstream queues $Q_1(t)$ and $Q_2(t)$ are also affected by links 4 and 5. The main result of this section is a criterion for exisitence priority vectors that stabilize the network:

\begin{Theorem}
\label{thm_2}
Consider a merge-diverge system. The merging junction has a priority vector $[\phi_1\ \phi_2]^T\in[0,1]^2$ satisfying \eqref{eq_phi} . The diverging junction has a discharging rule $\psi\in[0,1]$ satisfying \eqref{eq_psi}.
Then, there exists a non-empty set of static priorities $\phi\in[0,1]^2$ that stabilize the system if and only if
\begin{align}
\bar a_1<\min\{F_1,R_4\},\
\bar a_2<\min\{F_2,R_5\},\
\bar a_1+\bar a_2<F_3.
\label{eq_nominal2}
\end{align}
Furthermore, when \eqref{eq_nominal} holds, a set of stabilizing priority vectors is given by
\begin{align}
\Phi_2=\Big\{&\phi\in\Phi:\bar a_1-\min\{F_1,\phi_1F_3,R_4,(\phi_1/\phi_2)R_5\}<0,\nonumber\\ 
&\bar a_2-\min\{F_2,\phi_2F_3,R_5,(\phi_2/\phi_1)R_4\}<0\Big\}.
\label{eq_Phi2}
\end{align}
\end{Theorem}

Comparing the above result with Theorem~\ref{thm_1}, we can see that the diverge junction does not affect the existence of stabilizing priority vectors; however, it does affect the set of stabilizing priority vectors.
Apparently, $\Phi_2\subset\Phi_1$, where $\Phi_1$ is given by \eqref{eq_Phi1}.

\subsection{Proof of Theorem~\ref{thm_2}}

The necessity of \eqref{eq_nominal2} is apparent: if the fluid queuing process is stable, then the average inflow has to be less than the capacity for both traffic classes.

We next prove the sufficiency of \eqref{eq_nominal2}.
We only consider the case where $a_1^+>F_1$ and $a_2^+>F_2$; the other cases can be covered following analogous steps.
The proof is based on the following intermediate result:
\begin{Lemma}
\label{lmm_qtilde}
For the fluid queuing process over the merge-diverge system, suppose that $\phi\in\Phi_2$ and $a_k^+>F_k$, $k=1,2$. Then, for any $t>0$ and any $\tilde q_k\ge0$, $k=1,2$, if $Q_k(t)\ge\tilde q_k$, then 
\begin{align}
\psi_k(A(t),Q(t))\ge\tilde\psi_k^{\tilde q_k}:=\frac{\theta_k^{\tilde q_k}}{\Theta}
\label{eq_psik}
\end{align}
where $\theta_k^{\tilde q_k}$ is given by
\begin{align}
&\theta_k^{\tilde q_k}
=\nonumber\\
&\int_{s=0}^{\tilde q_k/(a_k^+-\min\{F_k,\phi_kF_3\})}\Big(\min\{F_k,\phi_kF_3\}-(\theta_s^k/\Theta)F_3\Big)ds.
\label{eq_theta_s}
\end{align}
\end{Lemma}

This result essentially states that under priority vector $\phi\in\Phi_2$, if there is a long queue in link 1 (resp. 2), then the fraction of class 1 (resp. 2) traffic in link 3 is subject to a lower bound.

\emph{Proof of Lemma~\ref{lmm_qtilde}:}
Suppose that $Q_1(t)=\theta_1\ge0$ and $Q_3^1(t)=\theta_3^1\ge0$ for a given $t>0$.
Since $\dot Q_1(s)\le a_1^+-\min\{F_1,\phi_1F_3\}$ for any $s$, we have
\begin{align*}
Q_1(s)\ge (a_1^+-\min\{F_1,\phi_1F_3\})(s-t_0)
\end{align*}
for each $s\in(t-\theta_1/(a_1^+-\min\{F_1,\phi_1F_3\}),t]$.

Thus, for each $s\in(t-\theta_1/(a_1^+-\min\{F_1,\phi_1F_3\}),t]$, $Q_1(s)>0$ and consequently
$
f_{12}(A(s),Q(s))=\min\{F_1,\phi_1F_3\}.
$
In addition, $\theta_s^1$ as specified in \eqref{eq_theta_s} is the solution to $Q_3^1(s)$ with the initial condition $Q_3^1(0)=0$, $Q_3^2(0)=\Theta$ and the constraints $Q_1(s)>0$ and $Q_2(s)$ for $t\in(0,\theta_1/(a_1^+-\min\{F_1,\phi_1F_3\})]$.
The above scenario is the one where $f_{34}$ would be minimized.
Hence, $f_{34}(A(t),Q(t))\ge f_{34}(\cdot,q)$, where $q=[q_1\ q_2\ q_3^1\ q_3^2]^T$ satisfies $q_1=\theta_1$ and $q_3^1=\theta^k_{\theta_1/(a_k^+-\min\{F_k,\phi_kF_3\})}$.

If $Q_3^1(t)=Q_3^2(t)=0$, then $\psi_1\ge\phi_1$, which naturally satisfies \eqref{eq_psik}. Otherwise,
$
\psi_1\ge\theta^k_{\theta_1/(a_k^+-\min\{F_k,\phi_kF_3\})}/\Theta
$.
This completes the proof.
\hfill$\blacksquare$

Lemma~\ref{lmm_qtilde} implies that if \eqref{eq_nominal2} holds, then the fluid queuing process over the merge-diverge system admits an invariant set $\mathcal M$ such that
\begin{align*}
\mathcal M=\cup_{\tilde q_1\ge0,\tilde q_2\ge0}\{q\in\mathcal Q:q_k\ge\tilde q_k,\psi_k(q)\ge\tilde\psi_k^{\tilde q_k},\ k=1,2\}.
\end{align*}
Lemma~\ref{lmm_qtilde} also implies that there exist $\hat q_1$ and $\hat q_2$ such that
\begin{align*}
&\forall q:q_1\ge\hat q_1,\ \psi_1(q)\ge 1-R_5/F_3,\\
&\forall q:q_2\ge\hat q_2,\ \psi_2(q)\ge 1-R_4/F_3,
\end{align*}
which enable us to derive Theorem~\ref{thm_2}:

\emph{Proof of Theorem~\ref{thm_2}:}
Let 
$$
x=\left[\begin{array}{c}
(q_1-\hat q_1)_++q_3^1\\
(q_2-\hat q_2)_++q_3^2
\end{array}\right]
$$
and consider the Lyapunov function
\begin{align*}
V_2(i,q):=&x^T
\left[\begin{array}{cc}
1 & \tilde\alpha\\
\tilde\alpha & \tilde\alpha^2
\end{array}\right]
x
+[\beta_i\ \tilde\alpha\beta_i]x,
\\
&\quad q\in\mathcal Q,\ i\in\{00,10,01,11\}
\end{align*}
where{
\begin{align*}
\tilde\alpha:=\frac12\Big(&\frac{\bar a_1}{\min\{F_2,\phi_2F_3,R_5,(\phi_2/\phi_1) R_4\}-\bar a_2}\\
&+\frac{\min\{F_1,\phi_1F_3,R_4,(\phi_1/\phi_2)R_5\}-\bar a_1}{\bar a_2}\Big)
\end{align*}}
and $\beta_i$ are specified by \eqref{eq_beta00}--\eqref{eq_beta11}.

Then, following procedures analogous to the proof of Proposition~\ref{prp_sufficient}, we can show that there exist $c>0$ and $d<\infty$ such that
$$\mathcal L_2V_2(a,q)\le-c|q|+d\quad\forall a\in\mathcal A,\ \forall q\in\mathcal M$$
which implies stability.
\hfill$\blacksquare$ 
\section{Numerical example}
\label{sec_simulate}

Now we use a numerical example to illustrate the results that we obtained in the previous sections.
Consider a merge-diverge network with the parameters in Table~\ref{tab_parameters}; the parameters make practical sense for the highway traffic setting.
\begin{table}[!ht]
  \caption{Model parameters}\label{tab_parameters}
  \begin{center}
    \begin{tabular}{|c|c|c|}
    \hline
    Parameter & Notation   &  Value                                   \\
    \hline
    Average inflow & $\bar a_1,\bar a_2$       &  $1200$ veh/hr                                 \\
    Merging link capacity & $F_1,F_2$    &  1500 veh/hr                           \\
    \begin{tabular}{c}Diverging link capacity \\(receiving flow)\end{tabular}
     & $R_4,R_5$     &  1400 veh/hr      \\
    \hline
    \end{tabular}
  \end{center}
\end{table}
In modern intelligent transportation systems, the Markovian inflow can be interpreted in a particular scenario, i.e. traffic flow with connected and autonomous vehicles traveling in platoons, or ``platooning'' \cite{al2010experimental}.
\begin{figure}[hbt]
\centering
\includegraphics[width=0.45\textwidth]{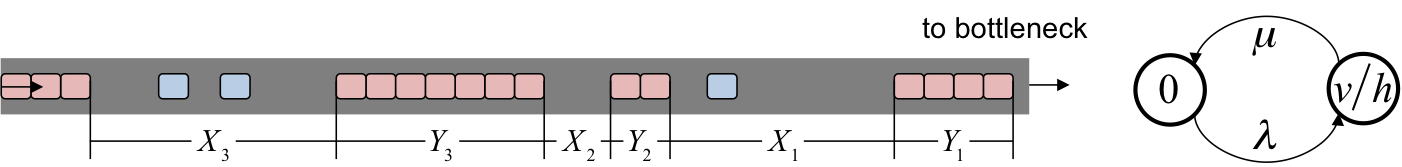}
\caption{Randomly arriving platoons cause Markovian switches in traffic flow.}
\label{fig_arrival}
\end{figure}
In this scenario, if we assume that (i) the headways between platoons are independent and identically distributed (IID) random variables $X$ with the cumulative distribution function (CDF)
$$\mathsf F_X(x)=1-e^{-\lambda x},\quad x\ge0,$$
and (ii) the lengths of platoons are IID with the CDF
$$\mathsf F_Y(y)=1-e^{-\lambda y},\quad y\ge0;$$
see Fig.~\ref{fig_arrival}. Furthermore, we assume that the background traffic flow is constant. Then, the inflow is a two-state Markov process \cite{jin2018modeling}.

We study the range of stabilizing priority vectors for various values of the capacity of the common link, $F_3$.
That is, for every given value of $F_3$ and given value of $\phi_1$, we check the stability conditions, i.e. verifying whether $\phi=[\phi_1\ \phi_2]^T$ is in the sets $\Phi_1'$, $\Phi_1$, and $\Phi_2$.

\begin{figure}[hbt]
\centering
\includegraphics[width=0.45\textwidth]{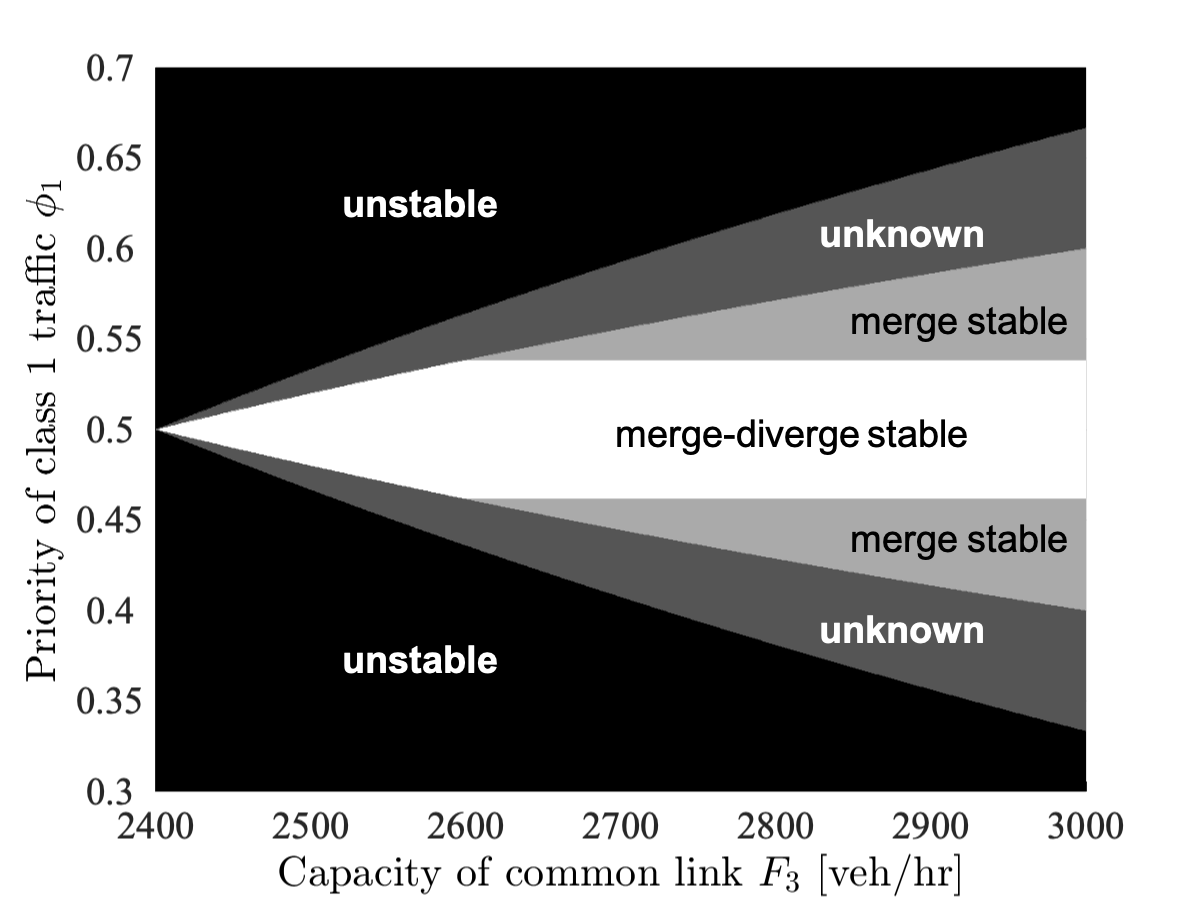}
\caption{Stability of the merge-diverge network under various priorities and various capacities of the common link.}
\label{stable}
\end{figure}

Fig.~\ref{stable} illustrates the results; the nomenclature is explained in Table~\ref{tab_stable}.
\begin{table}[!ht]
  \caption{Nomenclature for various regions in Fig.~\ref{stable}.}\label{tab_stable}
  \begin{center}
    \begin{tabular}{|c|c|c|c|}
    \hline
    Region & $\phi\in\Phi_0$?   &  $\phi\in\Phi_1$? & $\phi\in\Phi_2$?                                  \\
    \hline
    unstable & no       &  no & no                              \\
    unknown & yes    &  no & no                           \\
    merge stable & yes & yes & no \\
    merge-diverge stable & yes & yes & yes \\
    \hline
    \end{tabular}
  \end{center}
\end{table}
The following observations are noteworthy. First, there exists stabilizing priority vectors if and only if the common link has sufficient capacity to discharge both traffic classes, i.e. $F_3>2400=\bar a_1+\bar a_2$.
Second, there exist gaps (``unknown'') between the ``stable'' regions and the ``unstable'' region due to the gap between the necessary condition (characterized by $\Phi_0$) and the sufficient condition (characterized by $\Phi_1$) for stability.
Third, merge-diverge stability requires more restrictions on $\phi$ than merge stability alone.
Fourth, as long as $F_3$ is larger than a certain threshold (2600 in this example), the set of priority vectors stabilizing the merge-diverge network is insensitive to $F_3$; the reason is that in that range the downstream receiving flows $R_4$ and $R_5$ are the decisive quantities for stability. 
\section{Concluding remarks}
\label{sec_conclusion}

In this article, we studied the behavior of two traffic flows of distinct origins and destinations sharing a common link on their routes. Both flows are generated by a Markov process, and the delay is estimated using a fluid model.
We found that the way in which the limited space in the shared link is allocated to either traffic flow (characterized by the priority vector at the merging junction) plays a decisive role in the network's behavior.
In general, the fractional priority of a traffic flow should be in a neighborhood (which we quantitatively specify) of the inflow-to-capacity ratio of that flow.
Furthermore, although spillback may also happen at the diverging junction, it does not affect the stability or throughput or the network.

This work can serve as the basis for several directions of future work. 
First, modern traffic networks are equipped with real-time sensing and actuating capabilities. Therefore, the priority vector at the merging junction can be made responsive to real-time traffic condition. The advantage of a dynamic feedback priority vector is that it does not necessarily require accurate prediction of inflow or capacity.
Second, our analysis can be extended to more general networks, and approximated models can be developed for scalability.
Third, route choice model (road traffic) or routing algorithm (air traffic) can be added in the network extension as a second dimension of control capabilities. 

\bibliographystyle{IEEEtran}
\bibliography{bib_LJ}   
\end{document}